\documentclass{amsart}

\newcommand{\ignore}[1]{}

\usepackage{graphicx}

\newcommand{\dofigure}[1]{\begin{figure}#1 \end{figure}}
\newcommand{\bibtitle}[1]{\emph{#1}}
\newcommand{\figref}{Figure~\ref}
\newcommand{\dfn}[1]{\textbf{#1}}
\newcommand{\Erdos}{Erd\"os}

\newcommand{\nn}{\mathbb{N}}
\newcommand{\x}{ \times }
\newcommand{\st}{ : \,} % "such that" in set notation
\newcommand{\col}{\mathrm{col}}
\newcommand{\ch}{\mathrm{ch}}
\newcommand{\setsubtr}{\backslash}  % set subtraction

\newcommand{\knnn}{K_{n,n^n}}

\newtheorem{theorem}{Theorem}
\newcommand{\thm}[1]{ \begin{theorem} {#1} \end{theorem} }
\newtheorem{lemma}{Lemma}
\newcommand{\lem}[1]{ \begin{lemma} {#1} \end{lemma} }
\newtheorem{corollary}{Corollary}

\newtheorem{REMARK}{Remark}

\newtheorem{question}{Question}

\newcommand{\pf}[1]{\begin{proof}{#1}\end{proof}}

\newtheorem*{theoremx}{Theorem}
\newcommand{\thmx}[1]{ \begin{theoremx} {#1} \end{theoremx} }
\newtheorem*{lemmax}{Lemma}

\newtheorem*{corollaryx}{Corollary}
\newcommand{\corx}[1]{ \begin{corollaryx} {#1} \end{corollaryx} }
\newtheorem*{definitionx}{Definition}

\newtheorem{Question}{Question}

\begin{document}{

\title{List Coloring and $n$-Monophilic Graphs}
\author{Radoslav Kirov}
\address{
Department of Mathematics, 
University of Illinois at Urbana-Champaign,
Urbana, IL 61801, USA.
}
\author{Ramin Naimi}
\address{
Department of Mathematics, 
Occidental College,
Los Angeles, CA 90041, USA.
}
\date{28 April 2010}

\begin{abstract}
In 1990, Kostochka and Sidorenko
proposed studying the smallest number of list-colorings
of a graph $G$ among all assignments of lists of a given size $n$
to its vertices.
We say a graph $G$ is $n$-monophilic if
this number is minimized
when identical $n$-color lists
are assigned to all vertices of $G$.
Kostochka and Sidorenko observed that
all chordal graphs are $n$-monophilic for all $n$.
Donner (1992) showed that
every graph is $n$-monophilic
for all sufficiently large $n$.
We prove that all cycles are $n$-monophilic for all $n$;
we give a complete characterization of $2$-monophilic graphs
(which turns out to be similar to
the characterization of 2-choosable graphs given by
\Erdos, Rubin, and Taylor in 1980);
and for every $n$ we construct a graph
that is $n$-choosable but not $n$-monophilic.
\end{abstract}

\maketitle

\section{Introduction}

Suppose for each vertex $v$ of a graph $G$ we choose a list $L(v)$ of a fixed number $n$ of colors,
and then to each $v$ we assign a color chosen randomly from its color list $L(v)$.
If our goal is to maximize the probability of
getting the same color for at least two adjacent vertices,
then it seems intuitively plausible that
we should give every vertex of $G$ the same list.
But this turns out to be false for some graphs!
Graphs which do satisfy this property are called ``$n$-monophilic''
(defined more precisely below).
It is natural to ask: 
Which graphs are $n$-monophilic for a given  $n$?
This question has been open at least since 1990.

We work with finite, simple graphs,
and use the notation and terminology
of Diestel~\cite{Diestel}.
Given a graph $G = (V,E)$,
a \dfn{list assignment}
(resp.\ $n$-list assignment, $n \in \nn$)
for $G$ is a function
that assigns a subset (resp.\ $n$-subset) of $\nn$
to each vertex $v \in V$,
denoted $L(v)$.
Given a list assignment $L$ for $G$,
a (proper) \dfn{coloring} of $G$ from $L$
is a function $\gamma : G \to \nn$ such that
for each vertex $v \in V$,
$\gamma(v) \in L(v)$, and
for any pair of adjacent vertices $v$ and $w$,
$\gamma(v) \ne \gamma(w)$.
We denote the number of distinct colorings of $G$ from $L$ by
$\col(G,L)$.
In the special case where
$L(v) = [n]= \{1, \cdots, n \}$ for every $v \in V$,
we also write $\col(G,n)$ for $\col(G,L)$.
We say $G$ is \dfn{$n$-monophilic}
if
$\col(G,n) \le \col(G,L)$ for every $n$-list
assignment $L$ for $G$.
Clearly a graph is $n$-monophilic iff
each connected component of it is $n$-monophilic.
So we restrict attention to connected graphs only.

In 1990, Kostochka and Sidorenko \cite{KS}
proposed studying
the minimum value $f(n)$ attained by $\col(G,L)$
over all $n$-list assignments $L$ for a given graph $G$.
They observed that for chordal graphs
$f(n)$ equals the chromatic polynomial of $G$ evaluated at $n$;
i.e., chordal graphs are $n$-monophilic for all $n$.
In 1992 Donner \cite{Donner} showed that
for any fixed graph $G$,
$f(n)$ equals the chromatic polynomial of $G$
for all sufficiently large $n$;
i.e., every graph is $n$-monophilic for all sufficiently large $n$.
There appears to be no further literature on this subject since then.

A graph $G$ is said to be \dfn{$n$-colorable}
if $\col(G,n) \ge 1$;
and
$G$ is said to be \dfn{$n$-choosable}
(or \emph{$n$-list colorable})
if $\col(G,L) \ge 1$
for every $n$-list assignment $L$ for $G$.
The \dfn{chromatic number} of $G$, denoted $\chi(G)$,
is the smallest $n$ such that
$G$ is $n$-colorable.
The \dfn{list chromatic number} of $G$
(also called the \emph{choice number} of $G$),
denoted $\chi_l(G)$ (or $\ch(G)$),
is the smallest $n$ such that
$G$ is $n$-choosable.
Since $\chi$ and $\chi_l$ are well-known
and have been studied extensively,
it is interesting to compare
the concept of $n$-monophilic graphs to them.
The following are clear from definitions.
For every graph $G$,
\begin{enumerate}
\item
$\chi(G) \le \chi_l(G)$;
\item
if $n < \chi(G)$, then $G$ is $n$-monophilic;
\item
if $\chi(G) \le n < \chi_l(G)$, then $G$ is not $n$-monophilic.
\end{enumerate}
The interesting region is $\chi_l(G) \le n$,
which contains $n$-monophilic graphs
(e.g., all cycles and all chordal graphs),
as well as non-$n$-monophilic graphs (Section~\ref{Section:Choosable-but-not-monophilic}).

Deciding whether a graph is $n$-choosable turns out to be difficult.
Even deciding whether a given planar graph is $3$-choosable
is NP-hard \cite{Gutner}.
Thus one might expect the decision problem
for $n$-monophilic graphs to be NP-hard as well;
so a ``nice characterization'' 
(i.e., one that would lead to a polynomial time decision algorithm)
of $n$-monophilic graphs might not exist.
In this paper we prove that
all cycles are $n$-monophilic for all $n$,
and $G$ is not $2$-monophilic
iff all its cycles are even and 
it contains at least two cycles whose union is not $K_{2,3}$.
This characterization of 2-monophilic graphs 
is fairly similar to that given by \Erdos, Rubin, and Taylor~\cite{ERT}.
But, as we show in Section~\ref {Section:Choosable-but-not-monophilic},
for every $n \ge 2$
there is a graph that is
$n$-choosable but not $n$-monophilic.

\section{Chordal graphs are $n$-monophilic}

A graph is \dfn{chordal} if
every cycle in it of length greater than 3 has a chord.
Kostochka and Sidorenko~\cite{KS} observed that
all chordal graphs are $n$-monophilic for all $n$.
Because the proof is short,
we include it below.
Observe that if $H$ is a subgraph of $G$,
then $L$ restricts in a natural way to
give a list assignment for $H$,
and $\col(H,L)$ denotes
the number of colorings of $H$
from this restricted list assignment.

\lem{
\label{Lemma:Cone-on-clique}
Let $G$ be an $n$-monophilic graph,
and suppose $v_1, \cdots, v_k$
induce a complete subgraph of $G$.
Let $G'$ be the graph obtained from $G$
by adding a new vertex and connecting it to
$v_1, \cdots, v_k$.
Then $G'$ is $n$-monophilic.
}

\pf{
Let $L$ be an $n$-list assignment for $G'$.
If $n \le k$, then $\col(G',n) = 0$ and we are done.
So assume $n > k$.
Then each coloring of $G$ from $L$
extends to at least $n-k$ distinct colorings of $G'$ from $L$.
Hence
$\col(G',L)
\ge (n-k) \col(G,L)
\ge (n-k) \col(G,n)
= \col(G',n)$.
}

A graph has a \dfn{simplicial elimination ordering}
if its vertices can be ordered as
$v_1, \cdots, v_k$ such that
for each $v_i$
the subgraph induced by
$\{v_i\} \cup N(v_i) \cap \{v_1, \cdots, v_{i-1} \}$,
where $N(v_i)$ denotes the set of neighboring vertices of $v_i$,
is a complete graph.

\thmx{\emph{(Dirac \cite{Dirac})}
A graph is chordal iff it has a simplicial elimination ordering.
}

The above lemma and Dirac's Theorem
give us:

\corx{\emph{(Kostochka and Sidorenko \cite{KS})}
\label{Corollary:Chordal-graphs}
Every chordal graph is $n$-monophilic for every $n$.
}

Note that trees and complete graphs are chordal
and hence are $n$-monophilic for every $n$.

\section{Cycles are $n$-monophilic}
\label{Section:Cycles-are-monophilic}

In this section we show that
every $m$-cycle is $n$-monophilic for all $m, n$.
We first need some definitions.
Let $L$ be a list assignment for a graph $G$.
For $i = 1, \cdots, k$,
let $v_i$ be a vertex of $G$, and $c_i$ a color in $L(v_i)$.
Then $\col(G,L,v_1,c_1, \cdots, v_k,c_k)$
denotes the number of colorings of $G$ from $L$
which assign color $c_i$ to $v_i$,
$i = 1, \cdots, k$.
We say $L$ is \dfn{minimizing} for $G$
if $\col(G, L) \le \col(G, L')$
for every list assignment $L'$
where $|L'(v)|=|L(v)|$ for every $v$.

\lem{
\label{Lemma:Separating-edge}
Let $G_1$ and $G_2$ be disjoint subgraphs of a graph $G$,
with $v_i$ a vertex of $G_i$,
such that
$G= G_1 \cup G_2 + v_1 v_2$.
Let $L$ be a list assignment for $G$.
Then there exists a list assignment $L'$ such that 
$|L'(v)| = |L(v)|$ for every $v$,
$L'(v_1) \subseteq L'(v_2)$ or $L'(v_2) \subseteq L'(v_1)$,
and
$\col(G,L') \le \col(G,L)$.
Moreover,
the inequality is strict
provided there exist
$c_1 \in L(v_1)\setsubtr L(v_2)$
and
$c_2 \in L(v_2)\setsubtr L(v_1)$
with
$\col(G_1,L,v_1,c_1) \ne 0$
and
$\col(G_2,L,v_2,c_2) \ne 0$.}

\pf{
If $L(v_1) \subseteq L(v_2)$ or
$L(v_2) \subseteq L(v_1)$,
then there is nothing to show.
So we can assume there exist colors
$c_1 \in L(v_1) \setsubtr L(v_2)$ and
$c_2 \in L(v_2) \setsubtr L(v_1)$.

Let $L'$ be the list assignment that is identical to $L$
except that in the lists assigned to the vertices of $G_2$
every $c_1$ is replaced with $c_2$ and every $c_2$ with $c_1$.
Then, for each $c \ne c_1$ in $L(v_1)$,
$\col(G, L', v_1, c) = \col(G, L, v_1, c)$
(since $c \ne c_2$, as $c_2 \not\in L(v_1)$).
Furthermore,
\begin{equation}
\label{Equation:c1c2}
\begin{array}{rl}
\col(G, L', v_1, c_1) &
  =\col(G,L,v_1,c_1) - \col(G,L,v_1,c_1,v_2,c_2)  \\
 & = \col(G,L,v_1,c_1) - \col(G_1,L,v_1,c_1) \cdot \col(G_2,L,v_2,c_2)
\end{array}
\end{equation}
Hence, $\col(G, L', v_1, c_1) > \col(G,L,v_1,c_1)$
if $\col(G_1,L,v_1,c_1)$ and $\col(G_2,L,v_2,c_2)$ are both nonzero.

Now, by renaming $L'$ as $L$ and then
repeating this process as long as
$L(v_1) \not\subseteq L(v_2)$ and $L(v_2) \not\subseteq L(v_1)$,
we eventually obtain the desired $L'$.
}

The \dfn{length} of a path is the number of edges it contains.
For $n \ge 2$,
an \dfn{$(n,n-1)$-list assignment}
for a path of length at least one
is a function that assigns
$n$-color lists to the path's interior vertices, if any,
and $(n-1)$-color lists to its two terminal vertices.
Suppose the interior vertices of the path
have identical lists,
each of which contains as a subset
the $(n-1)$-color list
of each of the two terminal vertices.
If, in addition,
these two $(n-1)$-color lists are identical,
we say $L$ is \dfn{type~A},
and denote $\col(P, L)$ by $A_k$;
otherwise 
we say $L$ is \dfn{type~B},
and denote $\col(P, L)$ by $B_k$.
Note that, up to renaming colors,
all type A $(n,n-1)$-list assignments for a given path are equivalent,
and similarly for type B.

\lem{
\label{Lemma:Paths-are-monophilic}
Let $n \ge 2$, and let $L$ be an $(n,n-1)$-list assignment
for a path $P$ of length $k \ge 2$.
Then:
(a)~$A_k - B_k = (-1)^k$,
and
$A_k = \frac{n-1}{n}((n-1)^{k+1} + (-1)^k)$;
(b)~$\col(P,L) \ge \min(A_k, B_k)$;
and
(c)~for $n \ge 3$,
$L$ is minimizing only if
$k$ is odd and $L$ is type~A
or
$k$ is even and $L$ is type~B.
}

\begin{proof} 
Part (a):
Let $v$ be a terminal vertex of $P$,
and let $w$ be the vertex adjacent to $v$.
Suppose $L$ is type~A.
Then, for each color that we choose to assign to $v$,
there remains an $(n-1)$-color list of choices
for $w$, and this list is not the same as
the $(n-1)$-color list of the other terminal vertex of $P$.
Thus we get
\begin{equation}
\label{Equation:A-from-B}
A_k=(n-1)B_{k-1}
\end{equation}
By a similar (but slightly longer) reasoning, we see that
\begin{equation}
\label{Equation:B-from-AB}
B_k=A_{k-1} + (n-2)B_{k-1}
\end{equation}
Subtracting
(\ref{Equation:B-from-AB}) from (\ref{Equation:A-from-B})
gives
\begin{equation}
\label{Equation:A-minus-B-recursion}
A_k - B_k =(-1)(A_{k-1} - B_{k-1})
\end{equation}

Now, by direct calculation,
$A_1 = (n-1)(n-2)$,
and $B_1 = (n-2)^2 + (n-1)$.
It follows that $A_1 - B_1 = -1$,
which together with (\ref{Equation:A-minus-B-recursion})
inductively yield
\begin{equation}
\label{Equation:A-minus-B-closed}
A_k- B_k = (-1)^k
\end{equation}

Finally, combining
(\ref{Equation:A-from-B}) with
(\ref{Equation:A-minus-B-closed})
gives
$A_k = (n-1)(A_{k-1} + (-1)^k )$.
It follows by induction from the base case
$A_1 = (n-1)(n-2)$ that
$A_k = \frac{n-1}{n}((n-1)^{k+1} + (-1)^k)$.

Part (b):
This follows immediately from
Lemma~\ref{Lemma:Separating-edge}
and the definition of $A_k$ and $B_k$.

Part (c):
Assume $ k \ge 2$,
since the case $k=1$ is trivial.
Suppose, toward contradiction,
that $v_1, v_2$ are adjacent vertices of $P$ such that
$L(v_1) \not\subseteq L(v_2)$ and $L(v_2) \not\subseteq L(v_1)$.
Then, in the proof of Lemma~\ref{Lemma:Separating-edge},
the term
$\col(G_1,L,v_1,c_1) \cdot \col(G_2,L,v_2,c_2)$
being subtracted in equation~(\ref{Equation:c1c2})
is positive 
since each $G_i$ is now a path and $n \ge 3$.
This would imply that $L$ is not minimizing.
Hence, if $L$ is minimizing, it must be type A or type B.
The result now follows from equation~ (\ref{Equation:A-minus-B-closed}).

\end{proof}

Remark:
For $n =2$
a minimizing list need not be type A or type B;
examples are easy to construct.

\lem{
\label{Lemma:LargerListsGiveMoreColorings}
Let $L$ and $L'$ be distinct list assignments for a path $P$
such that for every vertex $v \in P$ we have
$L(v) \subseteq L'(v)$ and
$|L'(v)| \ge 2$.
Then $\col(P,L) < \col(P,L')$.
}

\begin{proof}
Since $L(v) \subseteq L'(v)$ for every $v$,
every coloring of $P$ from $L$ is
also a coloring of $P$ from $L'$.
And since $L$ and $L'$ are distinct,
for some vertex $w$ there is a color $c \in L'(w) \setsubtr L(w)$.
By hypothesis,$|L'(v)| \ge 2$ for every $v$;
hence  $\col(P,L',w,c) \ge 1$,
i.e., there is at least one coloring of $P$ from $L'$
that is not a coloring of $P$ from $L$.
The result follows.

\end{proof}

Let $L$ be a list assignment for a graph $G$,
and let $v_1, \cdots, v_k$ be vertices in $G$, where $k < |G|$.
Let $c_i \in L(v_i)$.
We define the list assignment $L_{c_1, \cdots, c_k}$
\dfn{induced} by $L$, $c_1, \cdots, c_k$
on the graph $H = G - \{v_1, \cdots, v_k\}$
by:
for every vertex $v \in H$,
$L_{c_1, \cdots, c_k}(v) = L(v) \setsubtr  \{c_i\ : v_i \in N(v) \}$
where $N(v)$ denotes the set of vertices in $G$ adjacent to $v$.
Then clearly
$\col(G,L,v_1,c_1, \cdots, v_k,c_k) = \col(H,L_{c_1, \cdots, c_k})$.

\thm{
\label{Theorem:cycles-are-monophilc}
Every cycle is $n$-monophilic for all $n \ge 2$.
}

\begin{proof} 
Let $C$ be a cycle of length $k \ge 3$.
Suppose we assign the color list $[n]$ to every vertex of $C$.
Then, by Lemma~\ref {Lemma:Paths-are-monophilic},
for each $c \in [n]$,
$\col(C,n,v,c) = A_{k-2}$.
Therefore $\col(C,n) = n A_{k-2}$.
Let $L$ be an $n$-list assignment
that does not assign identical lists to all vertices of $C$.
We will show $\col(C,L) \ge n A_{k-2}$.
Since $L$ does not assign identical lists to all vertices of $C$,
there are adjacent vertices $v$ and $w$ such that $L(v) \ne L(w)$.
Let $P$ be the path $C-v$.
We have two cases.

Case 1: $k$ is odd.
For each $c \in L(v)$, 
let $L_c$ be the list assignment induced on $P$,
and let $L'_c$ be an $(n,n-1)$-list assignment on $P$ obtained from $L_c$
by making, if necessary, the color lists of the endpoints of $P$ smaller.
Then, since $P$ has odd length $k-2$,
by Lemma~\ref {Lemma:Paths-are-monophilic},
for each $c \in L(v)$, 
$\col(P,L_c) \ge \col(P,L'_c) \ge A_{k-2}$.
Thus $\col(C,L) \ge n A_{k-2}$, as desired.

Case 2: $k$ is even.
First suppose $n \ge 3$ (we will treat the case $n=2$ separately).
If for every $c \in L(v)$ we have $\col(C,L,v,c) \ge A_{k-2}$,
then we are done.
So assume for some $c_0 \in L(v)$ we have $\col(C,L,v,{c_0}) < A_{k-2}$.
Then  $\col(P,L_{c_0}) < A_{k-2}$ since $\col(P,L_{c_0}) = \col(C,L,v,{c_0})$.
As $P$ has even length $k-2$,
it follows from
Lemma~\ref {Lemma:Paths-are-monophilic} part (c) and 
Lemma~\ref {Lemma:LargerListsGiveMoreColorings}
that 
$L_{c_0}$ must be a type B $(n,n-1)$-list assignment for $P$ and
$\col(P,L_{c_0}) = B_{k-2} = A_{k-2} -1$.
Let $u$ be the vertex adjacent to $v$ in $C-w$.
Then $c_0$ is in both $L_u$ and $L_w$,
and not in $L_{c_0}(u) \cup L_{c_0}(w) = L_{c_0}(x) = L(x)$,
where $x$ is any vertex in $P - \{u,w\}$.
Therefore, for each $c \ne c_0 \in L(v)$,
the induced list assignment $L_c$ on $P$ is not type B
because 
$c_0$ is in $L_c(u)$ and $L_c(w)$ but not in $L_c(x) = L(x)$.
Hence $\col(C,L,v,c) > B_{k-2}$, i.e.,
$\col(C,L,v,c) \ge A_{k-2}$.

Now, as $L(v) \ne L(w)$,
there exists an element $d \in L(v) \setsubtr L(w)$.
We show as follows that $\col(C,L,v,d) \ge A_{k-2}+1$.
Note that $L_d(w) = L(w)$ contains $n$ colors.
Let $L'_d$ be an $(n,n-1)$-list assignment for $P$
obtained from $L_d$ by removing one element other than $c_0$ from $L_d(w)$, 
and also one element from $L_d(u)$ if $|L(u)|=n$.
Since $c_0 \in L'_d(w)$, $L'_d$ is not a type B list assignment for $P$.
So, by Lemma~\ref{Lemma:Paths-are-monophilic},
$\col(P,L'_d) \ge B_{k-2}+ 1 = A_{k-2}$.
Hence, by Lemma~\ref{Lemma:LargerListsGiveMoreColorings},
$\col(C,L,v,d) = \col(P,L_d) > \col(P,L'_d) \ge A_{k-2}$,
as desired.

Thus we get

$$
\begin{array}{rl}
\col(C,L) & = \col(C,L,v,c_0) +  \col(C,L,v,d) +
      \sum\limits_{c \in L(v)\setsubtr \{c_0,d\}} \col(C,L,v,c) \\
 & \ge B_{k-2} + (A_{k-2} + 1)  + (n-2)A_{k-2}\\
 & = n A_{k-2} 
\end{array}
$$
as desired.

Now suppose $n=2$.
Clearly $\col(C,2) = 2$.
Let $Q$ be the path obtained by removing the edge $vw$ (but not its vertices) from $C$. 
We will show there are at least two colorings of $Q$ from $L$
that extend to colorings of $C$.
First, we need a definition. 
Let $x$ and $y$ be any two vertices in a graph $G$ 
with a given list assignment $M$.
We say that $c \in M(x)$ \dfn{forces} $d \in M(y)$
if $\col(G,M,x,c,y,d) \ge 1$ and
for every $d' \ne d \in M(y)$, $\col(G,M,x,c,y,d') = 0$.

Denote the vertices of $Q$ by $v_0, \cdots, v_k$,
where $v_0 = v$, $v_k =w$, and $v_i$ is adjacent to $v_{i+1}$ for $i = 0, \cdots, k-1$.
Suppose, toward contradiction, 
that each color in $L(v)$ forces a color in $L(w)$.
Then each color in $L(v_i)$ must be in $L(v_{i+1})$.
Hence $L(v_i) = L(v_{i+1})$.
But $L(v) \ne L(w)$.
So at least one of the colors in $L(v)$ forces no color in $L(w) $. 
Let $L(v) = \{\alpha, \beta\}$ and $L(w) = \{\gamma, \delta\}$. 
Then, without loss of generality, 
$\alpha$ forces neither $\gamma$ nor $\delta$. 
Therefore
$\col(Q,L,v,\alpha,w,\gamma)$
and 
$\col(Q,L,v,\alpha,w,\delta)$
are both nonzero,
since $\col(Q,L,v,\alpha) \ge 1$.
Now, if $\alpha$ is different from both $\gamma$ and $\delta$, 
then any coloring of $Q$ that assigns $\alpha$ to $v$
extends to a coloring of $C$, and we're done.
On the other hand, 
suppose $\alpha$ is not different from both $\gamma$ and $\delta$.
Then, without loss of generality, $\alpha = \gamma$. 
So any coloring of $Q$ with $\alpha$ assigned to $v$ and $\delta$ to $w$
extends to a coloring of $C$.
Also, $\beta \ne \gamma$ since $\beta \ne \alpha$. 
And $\beta \ne \delta$ since $\{\alpha, \beta\} \ne \{\gamma, \delta\}$. 
So any coloring of $Q$ with $\beta$ assigned to $v$
also extends to a coloring of $C$.
As $\col(Q,L,v,\beta) \ge 1$, 
we are done again.

\end{proof}

Note that although cycles are $2$-monophilic,
every even cycle has a minimizing $2$-list assignment
that does not assign the same list to every vertex:
assign the list $\{1,2\}$ to two adjacent vertices,
and the list $\{2,3\}$ to all the remaining vertices.

\section{A characterization of $2$-monophilic graphs}

The \dfn{core} of a connected graph $G$ is
the subgraph of $G$ obtained by
repeatedly deleting vertices of degree 1 until every remaining vertex has degree at least 2.

\lem{
\label{Lemma:monophilic-core}
A connected graph is $n$-monophilic iff its core is $n$-monophilic.
}

\begin{proof}
This is proved easily using
Lemma~\ref{Lemma:Separating-edge}
and induction on the number of vertices in the graph.

\end{proof}

Let $\theta_{a,b,c}$ denote the graph
consisting of two vertices
connected by three paths of lengths $a, b, c$
with mutually disjoint interiors.
In particular,
$\theta_{2,2,2}$ is the complete bipartite graph $K_{2,3}$.
In the paper by \Erdos, Rubin, and Taylor~\cite{ERT},
we find the following result by Rubin:

\thmx{
\emph{(A.\ L. Rubin)}
A connected graph is $2$-choosable iff its core is
a single vertex,
an even cycle,
or $\theta_{2,2,2m}$ for some $m \ge 1$.
}

We use this to prove that
a connected graph is $2$-monophilic iff
its core
is a single vertex,
is an even cycle,
is $K_{2,3}$,
or contains an odd cycle.

\lem{
\label{Lemma:K23-is-$2$-monophilic}
$K_{2,3}$ is $2$-monophilic.
}

\begin{proof}%\pf{
In \figref{Figure:K23-labeled}
the five vertices of $K_{2,3}$ have been labeled as $u,v,w,x,y$.
Since $K_{2,3}$ has no odd cycles,
$\col(K_{2,3},2) = 2$.
Let $L$ be a $2$-list assignment for $K_{2,3}$.
We will show that $\col(K_{2,3},L) \ge 2$.
We consider three cases,
depending on the number of colors
that $L(x)$ and $L(y)$ share.

\dofigure{[ht]
\centering
\includegraphics[width=30mm]{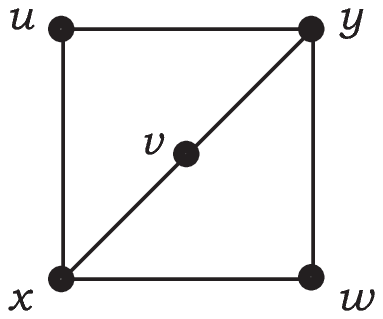}
\caption{$K_{2,3}$
\label{Figure:K23-labeled}}
}

\emph{Case 1.}
$|L(x) \cap L(y)| = 2$.
Then there are two ways to assign the same color to $x$ and $y$;
and for each way,
there is at least one way to color each of $u$, $v$, and $w$.
Hence $\col(K_{2,3},L) \ge 2$.

\emph{Case 2.}
$|L(x) \cap L(y)| = 1$.
Without loss of generality,
$L(x) = \{1,2\}$ and $L(y) = \{1,3\}$.
If at least one of the vertices $u,v,w$,
does not contain color~1 in its list,
then there are at least two distinct colorings of $K_{2,3}$
with color~1 assigned to both $x$ and $y$.
On the other hand,
if all three vertices $u,v,w$
contain color~1 in their lists,
then we can obtain one coloring
by assigning color~2 to $x$, 3 to $y$,
and 1 to $u,v,w$,
and another coloring
by assigning color~1 to both $x$ and $y$,
and using the second color in each of
the lists for $u,v,w$.

\emph{Case 3.}
$|L(x) \cap L(y)| = 0$.
Without loss of generality,
assume $L(x) = \{1,2\}$, $L(y) = \{3,4\}$.
Then there are four ways to color the pair $x,y$.
If at least two of these extend to a coloring of $K_{2,3}$,
we are done.
Otherwise, without loss of generality,
$L(u)=\{1,3\}$,
$L(v)=\{1,4\}$, and
$L(w)=\{2,3\}$.
Then
$(u,v,w,x,y)=(1,1,3,2,4)$ and
$(u,v,w,x,y)=(3,1,3,2,4)$
are two distinct colorings of $K_{2,3}$.

\end{proof}%}

\thm{
\label{Theorem:2-monophilc-characterization}
A connected graph is $2$-monophilic iff
its core
is a single vertex, is a cycle, is $K_{2,3}$,
or contains an odd cycle.
}

Equivalently:
A graph is \emph{not} $2$-monophilic
iff all its cycles are even
and it contains at least two cycles
whose union is not $K_{2,3}$.

\begin{proof} %\pf{

Clearly a single vertex
and a graph that contains an odd cycle
are both $2$-monophilic.
Also, by Theorem~\ref{Theorem:cycles-are-monophilc}
and Lemma~\ref{Lemma:K23-is-$2$-monophilic},
all cycles and $K_{2,3}$ are $2$-monophilic.
Using Lemma~\ref{Lemma:monophilic-core},
this gives us one direction of the theorem.

To prove the converse,
let $G$ be a $2$-monophilic graph.
If $G$ is not $2$-colorable,
then it must contain an odd cycle,
and we are done.
So assume $\chi(G)$ is  1 or 2.
Then $G$ is also $2$-choosable since it is $2$-monophilic.
So, by Rubin's theorem above,
it is enough to show that
for $m \ge 2$,
$\theta_{2,2,2m}$ is not $2$-monophilic.

\figref{Figure:nested-squares} shows
a $2$-list assignment $L$
for the case when $m = 2$.
When $m > 2$,
we add an even number of vertices
to the interior of the edge $uv$ in \figref{Figure:nested-squares}
and assign to each new vertex the list $\{1,2\}$.
It is then easy to check that for $m \ge 2$,
$\col(\theta_{2,2,2m},L) = 1 < 2 = \col(\theta_{2,2,2m},2)$,
as desired.

\dofigure{[ht]
\centering
\includegraphics[width=35mm]{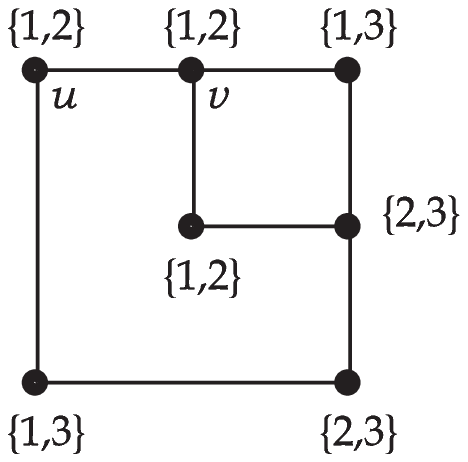}
\caption{$\col(\theta_{2,2,4},L)=1$.
\label{Figure:nested-squares}}
}

\end{proof} %} % end of proof of theorem

\section{Examples of $n$-choosable, non-$n$-monophilic graphs}
\label{Section:Choosable-but-not-monophilic}

Given the close similarity between Theorem~\ref {Theorem:2-monophilc-characterization}
and Rubin's theorem above,
it is natural to wonder how similar or different
the notions of $n$-choosalbe and $n$-monophilic are.
In this section,
for each $n \ge 2$
we construct a graph $H_{n}$
that is $n$-choosable but not $n$-monophilic.
To make the notation simpler,
we work with $H_{n+1}$ with $n \ge 1$
instead of $H_n$ with $n \ge 2$.

First, consider the complete bipartite graph $\knnn$.
Fix $n \ge 1$, and
denote the vertices of $\knnn$ by
$a_1, \cdots, a_n, b_1, \cdots, b_{n^n}$.
Let $L_0$ be an $n$-list assignment for $\knnn$ such that:
for all $i \ne j$,
$L_0(a_i) \cap L_0(a_j) = \emptyset$;
for all $k \ne l$,
$L_0(b_k) \ne L_0(b_l)$;
and
each $L_0(a_i)$ shares exactly one element with each $L_0(b_k)$.
Then there are $n^n$ distinct ways
to assign a color to each of $a_1, \cdots, a_n$,
and each of them will preclude
assigning a color to $b_k$ for some $k$.
It follows that
$\col(\knnn,L_0) = 0$.

Let $L'_0$ be an $n$-list assignment for $\knnn$
that is the same as $L_0$
except that its colors are renamed
so that colors $1, \cdots, n$
do not appear in any of its lists.
For each $j \in [n]$,
let $L_j$ be the $(n+1)$-list assignment for $\knnn$
given by $L_j(v) = L'_0(v) \cup \{j\}$
for every vertex $v \in \knnn$.
Let $x = \col(\knnn, L_j)$;
clearly $x$ is nonzero and independent of $j$.

Let $\{G_{i,j} \st i,j \in [n] \}$
be a set of $n^2$ disjoint copies of $\knnn$.
Let $p$ be the smallest integer such that
$n^p > x^{n^2}$.
Let $K_{n,p}$
be a complete bipartite graph with vertices
$v_1, \cdots, v_n, w_1, \cdots, w_p$.
We connect each $v_i$ to all vertices of
$G_{i,1}, \cdots, G_{i,n}$.
This describes the graph $H_{n+1}$.

\lem{
\label{Lemma:H-not-monophilic}
For all $n \ge 1$,
the graph $H_{n+1}$ is not $(n+1)$-monophilic.
}

\pf{
Define an $(n+1)$-list assignment $L$ for $H_{n+1}$ as follows.
For all $k \in [p]$,
$L(w_k) = \{n+1, n+2, \cdots, 2n+1 \}$;
for each $i \in [n]$,
$L(v_i) = [n] \cup \{ n+i\}$;
and on each $G_{i,j}$, $L = L_j$.

Let $\gamma$ be a coloring of $H_{n+1}$ from $L$.
Since $\col(\knnn,L_0) = 0$,
for each $i, j \in [n]$,
$\gamma$ must assign color $j$
to at least one vertex of $G_{i,j}$.
Hence for all $i \in [n]$, $\gamma(v_i) = n+i$;
and for all $k \in [p]$,
$\gamma(w_k) = 2n+1$.
It follows that
$\col(H_{n+1},L) = x^{n^2}$.

On the other hand,
$\col(H_{n+1}, n+1) \ge n^p$:
there are $n^p$ ways to color $w_1, \cdots, w_p$
from just $[n]$;
then assign color $n+1$ to every $v_i$;
and finally color every $G_{i,j}$
using colors 1 and 2.
Hence
$\col(H_{n+1},L) < \col(H_{n+1}, n+1)$,
as desired.

}

So it remains to show that $H_{n+1}$
is $(n+1)$-choosable.
We do this in the three following lemmas.
We say two list assignments $L$ and $L'$
for a graph $G$ are \dfn{equivalent}
if one can be obtained from the other by renaming colors and vertices,
i.e.,
there is a bijection $f : \nn \to \nn$
and an automorphism $\phi: G \to G$
such that for every vertex $v \in G$,
$L'(v) = f(L(\phi(v)))$.

\lem{
\label{Lemma:L0-is-unique}
Let $L$ be a list assignment for $\knnn$
such that for every vertex $v \in \knnn$,
$|L(v)| \ge n$.
If $\col(\knnn,L)=0$,
then $L$ is equivalent to $L_0$.
}

\pf{
Denote the two vertex-partitions of $\knnn$ by
$A$ and $B$,
with $|A|=n$ and $|B|=n^n$.
Suppose for some $a_1 \ne a_2$ in $A$,
$L(a_1) \cap L(a_2) \ne \emptyset$.
If we assign the same color to $a_1$ and $a_2$,
and to each $a' \ne a_1, a_2$
we assign a color to from $L(a')$,
then for every $b \in B$,
$L(b)$ contains at least one color that
was not assigned to any $a \in A$.
Hence $\col(\knnn, L) > 0$,
which is a contradiction.
So for all $a_1 \ne a_2$ in $A$,
$L(a_1) \cap L(a_2) = \emptyset$.

Now, suppose we have colored all $a \in A$.
Since $\knnn$ is not colorable from $L$,
there must exist some $b \in B$
whose color list is exactly
the $n$ colors we have chosen for the vertices in $A$.
Since there are only $n^n$ vertices in $B$,
if there were more than $n^n$ ways to color $A$,
then $\col(\knnn,L)$ would not be zero.
So each $L(a)$ must contain exactly $n$ colors,
and there are exactly $n^n$  ways to color $A$.
It follows that every $L(b)$
must contain exactly one color from $L(a)$ for each $a \in A$,
and no other colors;
Furthermore, distinct vertices in $B$ must have distinct lists.
This proves $L$ is equivalent to $L_0$.
}

\newcommand{\vconek}{K_{n,n^n,1}}

\lem{
\label{Lemma:Lj-is-unique}
Let $v$ denote the vertex in the one-element partition
of the complete tripartite graph $\vconek$.
Let $L$ be an $(n+1)$-list assignment for $\vconek$
such that $L(v) = [n+1]$.
Suppose for some $j \in [n]$,
$\col(\vconek,L,v,j) = 0$.
Then $L$ is equivalent to $L_j$;
and for all $i \ne j$,
$\col(\vconek,L,v,i) >0$.
}

\pf{
Let $L'$ be the list assignment for $\knnn$
obtained by deleting color $j$ from every list in
the restriction of $L$ to $\knnn$.
Then $\col(\knnn,L') = 0$.
So, by Lemma~\ref{Lemma:L0-is-unique},
$L'$ is equivalent to $L_0$,
and hence $L$ is equivalent to $L_j$.
If for some $i \ne j$, $\col(\vconek,L,v,i) = 0$,
then it would follow that
both $i$ and $j$ are in every list of $L$,
which contradicts the fact that
$L'$ is equivalent to $L_0$.
}

\lem{
For all $n \ge 1$,
$H_{n+1}$ is $(n+1)$-choosable.
}

\pf{
Let L be an $(n+1)$-list assignment for $H_{n+1}$.
For each $i,j \in [n]$,
let $G'_{i,j}$ be the subgraph of $H_{n+1}$
induced by $G_{i,j}$ and $v_i$.
By Lemma~\ref{Lemma:Lj-is-unique},
there is at most one color $c$ in $L(v_i)$
such that $\col(G'_{i,j},L,v_i,c)=0$.
Since $|L(v_i)|=n+1$,
there exists $c_i \in L(v_i)$ such that
for every $j \in [n]$,
$\col(G'_{i,j},L,v_i,c_i) \ne 0$.
Furthermore, since each $w_k$ has only $n$ neighbors,
$L(w_k) \setsubtr \{c_1, \cdots, c_n \}$ is non-empty.
Hence $\col(H_{n+1},L) \ne 0$.
}

\section{Questions}
\label{Section:Questions}

In this section we offer (and try to motivate)
two questions.
The Dinitz Conjecture,
proved by Galvin~\cite{Galvin},
states that
the line graph of the complete bipartite $K_{n,n}$ is $n$-choosable.
The List Coloring Conjecture
(which is open as of this writing),
generalizes the Dinitz Conjecture:
for every graph $G$,
$\chi_l(L(G)) = \chi(L(G))$,
where $L(G)$ denotes the line graph of $G$.

Note that the line graph of $K_{n,n}$
is isomorphic to the product $K_n \x K_n$,
where the product $G \x H$ is defined by
$V(G \x H) = V(G) \x V(H)$,
with two vertices $(g,h)$ and $(g',h')$
in $G \x H$ declared to be adjacent
if $g=g'$ and $h$ is adjacent to $h'$
or if $h=h'$ and $g$ is adjacent to $g'$.
Thus, another way to generalize the Dinitz Conjecture is:

\begin{Question}
Is the product of two $n$-monophilic graphs $n$-monophilic?
\end{Question}

For $n=2$ the answer to this question is \emph{No}:
Letting $P_i$ denote the path of length $i$,
it follows from Theorem~\ref{Theorem:2-monophilc-characterization}
that $P_2 \x P_3$ is not 2-monophilic,
while by Theorem~\ref{Lemma:Paths-are-monophilic}
every $P_i$ is 2-monophilic.
However, it is possible that the $n=2$ case is special
and for $n \ge 3$ the answer is \emph{Yes}.

\medskip

Donner's result \cite{Donner}
allows us to define the \dfn{monophilic number} of $G$,
denoted $\chi_m(G)$, in two possible natural ways:

\begin{enumerate}

\item
the smallest $n$ for which $G$ is $n$-colorable and $n$-monophilic, or

\item
the smallest $n$ such that $G$ is $n'$-monophilic for all $n' \ge
n$.

\end{enumerate}

We do not know whether or not these two definitions are equivalent;
it depends on the answer to the following:

\begin{Question}
If a graph is $n$-colorable and $n$-monophilic,
is it necessarily $(n+1)$-monophilic?
\end{Question}

\paragraph*{Acknowledgment}{
The authors thank Peter Keevash
of Queen Mary, University of London
for bringing the references
\cite{Donner} and \cite{KS}
to their attention.
The second named author is grateful to Caltech
for its hospitality while he worked there on this paper
during his sabbatical leave.
}

}
\begin{thebibliography}{1}


\bibitem{Diestel}
Reinhard Diestel.
\bibtitle{Graph Theory.}
Springer, Graduate Texts in Mathematics~\textbf{173} (1997).


\bibitem{Dirac}
G.\ A.\ Dirac.
\bibtitle{On rigid circuit graphs.}
Abhandlungen aus dem Mathematischen Seminar der Universitat Hamburg~\textbf{25} (1961), 71-76.

\bibitem{Donner}
Q.\ Donner.
\bibtitle{On the number of list-colorings.}
Journal of Graph Theorey~\textbf{16} (1992), 239-245.

\bibitem{ERT}
Paul \Erdos, Arthur L. Rubin, Herbert Taylor.
\bibtitle{Choosability in graphs.}
Congr.\ Numer.~\textbf{26} (1980), 125-157.

\bibitem{Galvin}
Fred Galvin.
\bibtitle{The List Chromatic Index of a Bipartite Multigraph.}
J.\ of Combinatorial Theory, Series B~\textbf{63} (1995), 153-158.


\bibitem{Gutner}
Shai Gutner.
\bibtitle{The complexity of planar graph choosability.}
Discrete Math.~\textbf{159} (1996) 119-130.


\bibitem{KS}
A.\ V.\ Kostochka, A.\ F.\ Sidorenko.
\bibtitle{Problem presented at the problem session.
Fourth Czechoslovak Symposium on Cominatorics,
Prachatice, 1990.}
Annals of Discrete Mathematics~\textbf{51} (1992), 380.


\bibitem{Tuza}
Z.\ Tuza.
\bibtitle{Graph colorings with local constraints --- A survery.}
Discussiones Mathematicae -- Graph Theory~\textbf{17}, No. 2 (1997), 161-228.



\end{thebibliography}
\end{document}